\documentclass[11pt]{article}
\usepackage{amstext,amssymb,amsmath,amsbsy}

\usepackage{hyperref}
\usepackage{amscd}
\usepackage{amsfonts}
\usepackage{indentfirst}
\usepackage{verbatim}
\usepackage{amsmath}
\usepackage{amsthm}
\usepackage{enumerate}
\usepackage{graphicx}
\usepackage{subfig}
\usepackage{color}
\usepackage[OT1]{fontenc}
\usepackage[latin1]{inputenc}
\usepackage[english]{babel}
\usepackage{amssymb}
\usepackage{dsfont}
\usepackage{esvect}

\newtheorem{theorem}{Theorem}
\newtheorem{lemma}{Lemma}

\newtheorem{proposition}{Proposition}

\newtheorem{remark}{Remark}

\numberwithin{equation}{section}

\newcommand{\proofend}{\hfill $\Box$ }

\newcommand{\dsp}{\displaystyle}

\newcommand{\e}{\varepsilon}

\newcommand{\supp}{\operatorname{supp}}

\newcommand{\dive}{\operatorname{div}}

\newcommand{\eps}{\varepsilon}

\newcommand{\mC}{\mathbb{C}}

\newcommand{\mN}{\mathbb{N}}

\newcommand{\mR}{\mathbb{R}}

\newcommand{\R}{{\cal R}}

\title{Superlensing using hyperbolic metamaterials: the scalar case}

\author{Eric Bonnetier \footnote{Universit\'e Grenoble-Alpes, Laboratoire Jean Kuntzmann, Grenoble, eric.bonnetier@imag.fr} \; \; 
Hoai-Minh Nguyen \footnote{EPFL SB MATHAA CAMA, Station 8,  CH-1015 Lausanne, hoai-minh.nguyen@epfl.ch}}

\begin{document}

\maketitle

\begin{abstract} This paper is devoted to superlensing using hyperbolic metamaterials: the possibility to image an arbitrary object using hyperbolic metamaterials without imposing any conditions on size of the object and the wave length.  To this end, two types of schemes are suggested and their analysis are given. The superlensing devices proposed are independent of the object. 
It is worth noting that the study of hyperbolic metamaterials is challenging due to the change type of modelling equations, elliptic in some regions, hyperbolic in some others. 
\end{abstract}

\section{Introduction}

Metamaterials are smart materials engineered to have properties that have 
not yet been found in nature.
They have recently attracted a lot of attention from the scientific community, 
not only because of potentially interesting applications,
but also because of challenges in understanding their peculiar properties. 

Negative index materials (NIMs) is an important class of such metamaterials.
Their study was initiated a few decades ago in the seminal paper of Veselago~\cite{Veselago},
in which he postulated the existence of such materials. New fabrication techniques now allow the construction
of NIMs at scales that are interesting for applications, and have made them a
very  active topic of investigation. One of the interesting properties of NIMs is superlensing, 
i.e., the possibility to beat  
the Rayleigh diffraction limit 
\footnote{The Rayleigh diffraction limit 
is on the resolution of lenses made of a standard {dielectric} material: 
the size of the smallest features in the images they produce is about  a half of the 
wavelength of the incident light. }: 
no constraint between the size of the object and the wavelength is imposed. 


Based on the theory of optical rays,
Veselago discovered that a slab lens of index~-1 could exhibit an unexpected 
superlensing property with no constraint on the size of the object
to be imaged~\cite{Veselago}.
Later studies by Nicorovici, McPhedran, and Milton~\cite{NicoroviciMcPhedranMilton94}, 
Pendry~\cite{PendryNegative, PendryCylindricalLenses}, 
Ramakrishna and Pendry in~\cite{PendryRamakrishna},
for constant isotropic objects and dipole sources,
showed similar properties for cylindrical lenses in the two dimensional  quasistatic regime, 
for the Veselago slab  and cylindrical lenses in the finite frequency regime, 
and for spherical lenses in the finite frequency regime.
Superlensing of arbitrary inhomogeneous objects using NIMs in the acoustic  and electromagnetic settings  
was established in~\cite{Ng-Superlensing, Ng-Superlensing-Maxwell} for related 
lens designs.  Other interesting properties of NIMs include cloaking using 
complementary media~\cite{LaiChenZhangChanComplementary, Ng-Negative-Cloaking, MinhLoc2}, 
cloaking a source via anomalous localized 
resonance~\cite{AmmariCiraoloKangLeeMilton, BouchitteSchweizer10, KohnLu, MiltonNicorovici, Ng-CALR, Ng-CALR-frequency, MinhLoc1}, 
and cloaking an arbitrary object via anomalous localized resonance~\cite{Ng-CALR-object}.  


In this paper, we are concerned with another type of metamaterials: hyperbolic metamaterials (HMMs).
These materials have quite promising potential applications to subwavelength imaging and focusing;
see \cite{Poddubny13} for a recent interesting survey on hyperbolic materials and their applications. 
We focus here on their superlensing properties. 
{The peculiar properties and the difficulties} in the study of NIMs come from (can be explained by) 
the fact that the equations modelling their behaviors have sign changing coefficients.
{
In contrast, the modeling of HHMs involve equations of changing type, elliptic in some regions, hyperbolic in others.}

We first describe a general setting concerning HMMs and point out some of their general properties.
Consider a standard medium that occupies a region $\Omega$ of $\mR^d$ ($d=2, 3$) 
with material constant $A$,
except  for a subset $D$ in which the material is hyperbolic  
with material constant $A^H$ in  the quasistatic regime
(the finite frequency regime is also considered in this paper and is discussed later).  
Thus, $A^H$ is a symmetric hyperbolic matrix-valued function defined in $D$  
and $A$ is a symmetric uniformly elliptic matrix-valued function defined in $\Omega \setminus D$.  
Since metamaterials usually contain damping (metallic) elements,
it is also relevant to assume that the medium in $D$ is lossy 
(some of its electromagnetic energy is dissipated as heat) and study the situation as the loss goes to 0. 
The loss can be taken into account by adding an imaginary
part of amplitude $\delta > 0$ to $A^H$.
With the loss, the medium in the whole of $\Omega$ is 
thus characterized by the matrix-valued function $A_\delta$ 
defined by 
\begin{equation}
A_\delta = \left\{\begin{array}{cl}
A &  \mbox{ in } \Omega \setminus D, \\[6pt]
A^H - i \delta I  & \mbox{ in } D.  
\end{array}\right.
\end{equation}
For a given (source) function $f \in L^2(\Omega)$, the propagation of light/sound
is modeled in the quasistatic regime by the equation
\begin{equation}\label{H-E}
\dive(A_\delta \nabla u_\delta ) = f  \mbox{ in } \Omega, 
\end{equation}
with an appropriate boundary condition on $\partial \Omega$. 

Understanding the behaviour of $u_\delta$ as $\delta \to 0_+$ is 
a difficult question in general due to two facts. 
Firstly, equation \eqref{H-E} has both elliptic (in $\Omega \setminus D$) and hyperbolic (in $D$) characters. 
It is hence out of the range of the standard theory of elliptic and hyperbolic equations. 
Secondly, even if \eqref{H-E} is of hyperbolic character in $D$, the situation is far from standard 
since the problem in $D$ is not an initial boundary problem. 
There are constraints on both the Dirichlet and Neumann boundary conditions (the transmission conditions).
As a consequence, equation \eqref{H-E} is very unstable (see Section~\ref{sect-toy1} for a concrete example).

\medskip

In this paper, we  study superlensing using HMMs. 
The use of hyperbolic media in the construction of lenses was  suggested by Jacob et al. in~\cite{Jacob06} 
and was experimentally verified by Liu et al. in ~\cite{Liu07}. 
The proposal of~\cite{Jacob06} concerns cylindrical lenses in which the hyperbolic material is given 
in standard polar coordinates by  
\begin{equation}\label{Jacob}
A^H = a_\theta e_\theta \times e_\theta - a_r e_r \times e_r,  
\end{equation}
where $a_\theta$ and $a_r$ are positive constants \footnote{It seems to us that in their proposal these constants can be chosen quite freely.}. Denoting the inner radius and the outer radius of the cylinder respectively by $r_1$  and $r_2$, 
Jacob et al. argued that 
\begin{equation}\label{Jacob-resolution}
\mbox{ the resolution is } \frac{r_1}{r_2} \lambda,  
\end{equation}
where $\lambda$ is the wave number. They supported their prediction by numerical simulations. 

The goal of our paper is to go beyond the resolution problem to achieve superlensing using HMMs 
as discussed in~\cite{Ng-Superlensing, Ng-Superlensing-Maxwell} in the context of NIMs, i.e., 
{to be able to image an object without imposing restrictions on the ratio between
its size and the wavelength of the incident light.}
We propose two constructions for superlensing, which are based on two  different mechanisms, 
inspired by two basic properties of {the} one dimensional wave {equation}.

The first mechanism is based on the following simple observation. 
Let $u$ be a smooth solution of the system
\begin{equation}\label{wave1-ole}
\left\{\begin{array}{c}
\partial_{tt}^2 u (t, x) - \partial_{xx}^2 u (t, x) = 0  \mbox{ in } \mR_+ \times [0, 2 \pi], \\[6pt]
u(t, \cdot) \mbox{ is $2\pi$-periodic}. 
\end{array}\right. 
\end{equation}
Then $u$ can be written {in} the form 
\begin{equation*}
u(t, x) = \sum_{-\infty}^\infty \sum_{\pm} a_{n, \pm} e^{in t \pm n x} \mbox{ in } \mR_+ \times [0, 2 \pi],  
\end{equation*}
for some  constant $a_{n, \pm} \in \mC$. This implies  
\begin{equation}\label{P-wave1}
u(t, \cdot) = u(t + 2 \pi, \cdot) \mbox{ for all } t \ge 0. 
\end{equation}
The key point here is that \eqref{P-wave1} holds for arbitrary Cauchy data at $t = 0$. 
Based on this observation, we propose the following two dimensional superlensing device in the annulus 
$B_{r_2} \setminus B_{r_1}$:
\begin{equation}\label{scheme1}
A^H =\frac{1}{r} e_r \times e_r - r e_\theta \times e_\theta \mbox{ in } B_{r_2} \setminus B_{r_1}, 
\end{equation}
under the requirement that 
\begin{equation}\label{cond-r1-r2}
r_2 - r_1 \in 2 \pi \mN_+
\end{equation}
(see \eqref{scheme1-k-2} for a three dimensional scheme in the finite frequency regime which is related to this observation). 
Here and in what follows $B_r$ denotes the open ball in $\mR^d$ centered at the origin 
and of radius $r$. We also use the standard notations for polar coordinates in two dimensions and spherical coordinates in three dimensions hereafter.
Given the form~\eqref{scheme1} of $A^H$, one can verify that 
$$
\dive(A^H \nabla u) = \frac{1}{r} (\partial_{rr}^2 u - \partial_{\theta \theta}^2 u) \mbox{ in } B_{r_2} \setminus B_{r_1}. 
$$
Hence, if $u$ is a  solution to the equation  $\dive(A^H \nabla u) = 0 $ in $B_{r_2} \setminus B_{r_1}$ then 
\begin{equation}\label{wave1-ole}
 \partial_{rr}^2 u - \partial_{\theta \theta}^2 u  = 0 \mbox{ in } B_{r_2} \setminus B_{r_1}. 
\end{equation}
It follows from \eqref{cond-r1-r2} that
\begin{equation}\label{gluing-cond} 
u(r_2 x/|x|) = u(r_1 x/ |x|) \quad \mbox{ and } \quad \partial_r u(r_2 x/|x|) = \partial_r u(r_1 x/ |x|).
\end{equation} 
This in turn implies the magnification of the medium contained
inside $B_{r_1}$ by a factor $r_2/r_1$ (the precise meaning is given in Theorem~\ref{thm-main}). Inspired by \eqref{P-wave1}, we call this scheme ``tuned superlensing" using HMMs. 

\medskip 
Our second class of superlensing devices is inspired by another 
observation concerning the one dimensional  wave equation.
Given $T> 0$, let $u$ be a solution with appropriate regularity to the system
\begin{equation}\label{wave2}
\left\{\begin{array}{c}
\partial_{tt}^2 u - \partial_{xx}^2 u = 0  \mbox{ in } (-T, 0) \times [0, 2 \pi], \\[6pt]
- \partial_{tt}^2 u  + \partial_{xx}^2 u = 0  \mbox{ in } (0, T) \times [0, 2 \pi], \\[6pt]
u \mbox{ is $2\pi$-periodic w.r.t. $x$}, \\[6pt]
u (0_+, \cdot) = u(0_-, \cdot ), \,  \partial_t u (0_+, \cdot) = - \partial_t u(0_-, \cdot )   \mbox{ in }  [0, 2 \pi]. 
\end{array}\right. 
\end{equation}
Then 
\begin{equation}\label{P-wave2}
u(t, x) = u(-t, x) \mbox{ for } (t, x) \in (0,T) \times [0, 2 \pi].  
\end{equation}
Indeed, set 
\begin{equation*}
v(t, x) = u(-t, x) \quad \mbox{ and } \quad w(t, x) = v(t, x) - u(t, x)  \mbox{ for } (t, x) \in (0, T) \times (0, 2 \pi). 
\end{equation*}
Then 
\begin{equation*}
\left\{\begin{array}{c}
\partial_{tt}^2 w  - \partial_{xx}^2 w = 0  \mbox{ in } (0, T) \times [0, 2 \pi], \\[6pt]
w(\cdot, 0) = w(\cdot, 2 \pi) = 0 \mbox{ in } (0, T), \\[6pt]
w \mbox{ is $2\pi$-periodic w.r.t. $x$}, \\[6pt]
w (0_+, \cdot) =  \partial_t w (0_+, \cdot) = 0   \mbox{ in }  [0, 2 \pi]. 
\end{array}\right. 
\end{equation*}
Therefore, $w = 0$ in $(0, T) \times (0, 2 \pi)$ by the uniqueness of the Cauchy problem for the wave equation;  which implies that $u(t, x) = u(-t, x)$ for $(t, x) \in (0, T) \times (0, 2 \pi)$ as mentioned. 
Based on this observation, we propose the following  superlensing device in $B_{r_2} \setminus B_{r_1}$ in both two and three dimensions, with $r_m = (r_1 + r_2)/ 2$: 
\begin{equation*}
A^H= \left\{ \begin{array}{cl} \dsp \frac{1}{r}e_r \otimes e_r - r e_\theta \otimes e_\theta & \mbox{ in } B_{r_2} \setminus B_{r_m}, \\[6pt]
\dsp- \frac{1}{r} e_r \otimes e_r + r e_\theta \otimes e_\theta  & \mbox{ in }  B_{r_m} \setminus B_{r_1}, 
\end{array} \right. \quad \mbox{ for } d = 2
\end{equation*}
and 
\begin{equation*}
A^H= \left\{ \begin{array}{cl} \dsp \frac{1}{r^2} e_r \otimes e_r - ( e_\theta \otimes e_\theta + e_\varphi \otimes e_\varphi)  & \mbox{ in } B_{r_2} \setminus B_{r_m}, \\[6pt]
-  \dsp \frac{1}{r^2}e_r \otimes e_r + ( e_\theta \otimes e_\theta + e_\varphi \otimes e_\varphi)   & \mbox{ in }  B_{r_m} \setminus B_{r_1}, 
\end{array} \right.  \quad \mbox{ for } d = 3. 
\end{equation*}
In a compact form, one has 
\begin{equation}\label{scheme2}
A^H= \left\{ \begin{array}{cl} \dsp \frac{1}{r^{d-1}}e_r \otimes e_r - r^{3-d} (I - e_r \otimes 
e_r)  & \mbox{ in } B_{r_2} \setminus B_{r_m}, \\[6pt]
\dsp - \frac{1}{r^{d-1}} e_r \otimes e_r + r^{3-d} (I - e_r \otimes e_r)  & \mbox{ in }  B_{r_m} \setminus B_{r_1}. 
\end{array} \right.
\end{equation}
From the definition of $A^H$ in \eqref{scheme2}, we have 
$$
\dive (A^H \nabla u) = \frac{1}{r^{d-1}} \Big( \partial^2_{rr} u - \Delta_{\partial B_1} u \Big) \mbox{ in } B_{r_2} \setminus B_{r_m},  
$$
and 
$$
\dive (A^H \nabla u) = -  \frac{1}{r^{d-1}} \Big(  \partial^2_{rr} u - \Delta_{\partial B_1} u \Big) \mbox{ in } B_{r_m} \setminus B_{r_1},  
$$
where $\Delta_{\partial B_1}$ denotes the Laplace-Beltrami operator on the unit sphere of $\mR^d$.
Hence, if $u$ is an appropriate solution to the equation 
$\dive (A^H \nabla u) = 0$ in $B_{r_2} \setminus B_{r_1}$, 
then, by taking into account the transmission conditions on $\partial B_{r_m}$, one has  
\begin{equation}\label{scheme2-11}
\left\{ \begin{array}{cl} \dsp  \partial^2_{rr} u  - \Delta_{\partial B_1} u = 0  & \mbox{ in } B_{r_2} \setminus B_{r_m}, \\[6pt]
\dsp -  \partial^2_{rr} u +  \Delta_{\partial B_1} u = 0 & \mbox{ in }  B_{r_m} \setminus B_{r_1}, \\[6pt]
\dsp u \big|_{B_{r_2} \setminus B_{r_m}} = u \big|_{B_{r_m} \setminus B_{r_1}}, \; \; \partial_r u \big|_{B_{r_2} \setminus B_{r_m}} = - \partial_r u \big|_{B_{r_m} \setminus B_{r_1}} &  \mbox{ on } \partial B_{r_m}. 
\end{array} \right.
\end{equation}
As in \eqref{P-wave2}, one derives that 
$$
u\big((s + r_m) \hat x \big) = u\big((r_m -s) \hat x \big) \mbox{ for } \hat x \in \partial B_{1}, s \in (0, r_2 - r_m);  
$$
which yields 
\begin{equation}\label{gluing-cond} 
u(r_2 \hat x) = u(r_1 \hat x) \quad \mbox{ and } \quad \partial_r u(r_2 \hat x) = - \partial_r u(r_1 \hat x) \mbox{ for } \hat x \in \partial B_1.
\end{equation}  
This in turn implies the magnification of the medium contained
inside $B_{r_1}$ by a factor $r_2/r_1$ (the precise meaning is given in Theorem~\ref{thm-main}).
In contrast with the first proposal \eqref{scheme1} where \eqref{cond-r1-r2} is required, 
we do not impose any conditions on $r_1$ and $r_2$ for the second scheme \eqref{scheme2}.  
We call  this method ``superlensing using HHMs via complementary property".  
The idea of using reflection takes roots in the work 
of the second author~\cite{Ng-Complementary}. Similar ideas were used in the study properties of NIMs 
such as superlensing~\cite{Ng-Superlensing, Ng-Superlensing-Maxwell}, 
cloaking~\cite{Ng-Negative-Cloaking, MinhLoc2}, 
cloaking  via anomalous localized resonance in~\cite{Ng-CALR, Ng-CALR-frequency, Ng-CALR-object, MinhLoc1}, 
and the stability of NIMs in~\cite{Ng-WP}. 
Nevertheless, the superlensing properties of NIMs and HMMs are based on two different phenomena: 
the unique continuation principle for NIMs,
and the uniqueness of the Cauchy problem for the wave equation for HMMs. 

\medskip
Suppose that an object to-be-magnified, located in $B_{r_1}$, is characterized by a symmetric uniformly elliptic   
matrix-valued function $a$. Throughout the paper, we assume that \footnote{This assumption is used to obtain enough regularity for solutions to deal with wave equations.}
\begin{equation}\label{smoothness}
a \mbox{ is of class $C^1$ in a neighborhood of $\partial B_{r_1}$.}
\end{equation}
Suppose that outside $B_{r_2}$ the medium is homogeneous. 
The whole system (taking loss into account) is then given by
\begin{equation}\label{medium-introduction}
A_\delta = \left\{\begin{array}{cl} I & \mbox{ in } \Omega \setminus B_{R_2},\\[6pt]
A^{H} - i \delta I &  \mbox{ in } B_{r_2} \setminus B_{r_1}, \\[6pt]
a & \mbox{ in } B_{r_1},  
\end{array}\right.
\end{equation}
where $A^H$ is defined either {by} \eqref{scheme1}-\eqref{cond-r1-r2} 
if $d=2$ or {by} \eqref{scheme2} if $d=2, 3$. 
Set 
\begin{equation}\label{def-H}
H^1_{m} (\Omega): = \Big\{u \in H^1(\Omega); \int_{\partial \Omega} u = 0 \Big\}. 
\end{equation}
One of the main results of this paper{,} 
stated here in the quasistatic regime{,} is 

\begin{theorem} \label{thm-main} 
Let $d=2, 3$, $0 < \delta < 1$, $0< r_1 < r_2$,  $\Omega$ be a smooth bounded connected open subset of $\mR^d$, and  let $f \in L^2(\Omega)$ with $\dsp \int_{\Omega} f = 0$. 
Assume that $B_{r_2} \subset \subset \Omega$ and  $\supp f \subset \Omega \setminus B_{r_2}$. Let $u_\delta \in H^1_{m}(\Omega)$ be the unique solution to the system 
\begin{equation}\label{eq-u-delta}
\left\{\begin{array}{cl}
\dive(A_\delta \nabla u_\delta) = f &  \mbox{ in } \Omega\\[6pt]
\partial_\nu u_\delta = 0 &  \mbox{ on } \partial \Omega, 
\end{array}\right.
\end{equation}
where $A_\delta$ is given by \eqref{medium-introduction}. 
We have 
\begin{equation}\label{thm-main-statement1}
\|u_\delta  \|_{H^1(\Omega)} \le C  \| f\|_{L^2(\Omega)} \quad \mbox{ and } \quad 
u_\delta \to u_0 \;\textrm{strongly in}\; H^1(\Omega), 
\end{equation}
where $u_0 \in H^1_{m}(\Omega)$ is 
the unique solution to~\eqref{eq-u-delta} with $\delta =0${,} 
and $C$ is a positive constant independent of $f$ and $\delta$. 
Moreover, $u_0 = \hat u$ in $\Omega \setminus B_{r_2}$ where $\hat u \in H^1_{m}(\Omega)$ 
is the unique solution to the system 
\begin{equation}\label{def-hA-hu}
\left\{\begin{array}{cl}
\dive(\hat A  \nabla \hat u) = f &  \mbox{ in } \Omega\\[6pt]
\partial_\nu \hat u = 0 &  \mbox{ on } \partial \Omega{,}
\end{array}\right. \quad
\mbox{ where }  \quad
\hat A(x) = \left\{\begin{array}{cl}  I  & \mbox{ in } \Omega \setminus B_{r_2},  \\[6pt]
\dsp \frac{r_1^{d-2}}{r_2^{d-2}} a\Big( \frac{r_1 }{ r_2 }  x\Big)& \mbox{ in } B_{r_2}.
\end{array}\right.
\end{equation}
\end{theorem}
\medskip

The well-posedness and the stability of \eqref{eq-u-delta} are established in Lemma \ref{lem1}.  
The existence and uniqueness of $u_0$ are a part of Theorem~\ref{thm-main}.  
Since $f$ is arbitrary with support in $\Omega \setminus B_{r_2}$, it follows from the definition 
of $\hat A$ that the object in $B_{r_1}$ is magnified by a factor $r_2/ r_1$.  
It is worth noting that $a$ can be an arbitrary function inside $B_{r_1}$,
provided it is uniformly elliptic and smooth near $\partial B_{r_1}$. We emphasize here  that the lens is independent of the object. 
 
\medskip

The paper is organized as follows. Section 2 is devoted to tuned superlensing
via HMMs. There, besides the proof of Theorem~\ref{thm-main}, 
where $A^H$ is given by \eqref{scheme1}-\eqref{cond-r1-r2},  
we also discuss a two dimensional  variant in the finite frequency domain (Theorem~\ref{thm-scheme1-k}),
and a result for the three dimensional  finite frequency regime, 
where $A^H$ is strictly hyperbolic (Theorem~\ref{thm-scheme1-k-2}).  
Section~\ref{sect-complementary} concerns superlensing using HMMs via {the} 
complementary property. 
In this section, we consider coefficients $A^H$ given by~\eqref{scheme2},
and prove a finite frequency generalization of Theorem~\ref{thm-main} 
(Theorem~\ref{thm-scheme2-k}).  
Finally,  in Section 4, we construct HMMs  with the required properties, 
as {limits as $\delta \to 0$ of} effective media obtained from the homogenization
of composite structures, mixtures of a dielectric and a ``real metal".
{
Numerical simulations of some of the results presented in our paper are presented in \cite{DHN}.}
\medskip 



\section{Tuned superlensing using HMMs}\label{sect-tuned}

In this section, we first present two lemmas on the stability of \eqref{H-E} and \eqref{def-hA-hu} and their variants in 
the finite frequency regime. 
In the second part, we discuss a toy model which illustrates tuned superlensing with hyperbolic media. 
Finally, we give a proof of Theorem~\ref{thm-main} when $A^H$ is given by~\eqref{scheme1}-\eqref{cond-r1-r2}, 
and we discuss its variants in the finite frequency case. 


\subsection{Two useful lemmas}

We first establish the following lemma which implies the well-posedness  of \eqref{H-E}. In what follows, for a subset $D$ of $\mR^d$, $\mathds{1}_{D}$ denotes its characteristic function. We have 

\begin{lemma} \label{lem1} Let $d  =2 , 3$, $k \ge  0$, $\delta_0>0$,  $0< \delta < \delta_0$.
Let $D \subset \subset \Omega$ be two  smooth bounded connected open subsets of $\mR^d$. Let $A$ be a bounded matrix-valued function defined in $\Omega$ such that $A$ is uniformly elliptic in $\Omega \setminus D$, $A$ is piecewise $C^1$ in $\Omega$,  and let $\Sigma$ be a complex bounded function such that $\Im(\Sigma) \ge 0$. Set  
\begin{equation}\label{def_coeff}
A_\delta(x) = A(x) - i \delta  \mathds{1}_{D}(x) I \mbox { and } \Sigma_\delta(x) = \Sigma(x) + i \delta \mathds{1}_{D}(x)\mbox{ in } \Omega. 
\end{equation}
Let $g_\delta \in  [H^1(\Omega)]^*$, the {dual space} of $H^1(\Omega)$,  
and in the case $k = 0$, assume in addition that $\int_{\Omega} g_\delta = 0$. 
There exists a unique solution $v_\delta \in H^1(\Omega)$ if $k > 0$ 
(respectively $v_\delta \in H^1_{m}(\Omega)$ if $k = 0$) to the system 
\begin{equation}\label{lem1-sys}
\left\{\begin{array}{cl}
\dive(A_\delta \nabla v_\delta) + k^2 \Sigma_\delta v_\delta =g_\delta &  \mbox{ in } \Omega, \\[6pt]
A \nabla v_\delta \cdot \nu - i k v_\delta  = 0 & \mbox{ on } \partial \Omega. 
\end{array}\right.
\end{equation}
Moreover, 
\begin{equation}\label{lem1-part1}
\| v_\delta \|_{H^1(\Omega)}^2 \le \frac{C}{\delta} \left| \int_{\Omega} g_\delta \bar v_\delta \right| + \|g_\delta \|_{[H^1(\Omega)]^*}^2, 
\end{equation}
for some positive constant $C$ depending only on $\Omega$, $D$, and $k$. Consequently, 
\begin{equation}\label{lem1-part2}
\| v_\delta \|_{H^1(\Omega)} \le \frac{C}{\delta} \| g_\delta\|_{[H^1(\Omega)]^*}. 
\end{equation}
\end{lemma} 

\noindent{\bf Proof.}  We only prove the result for $k>0$. The case $k=0$ follows similarly and 
{is} left to the reader. 
The proof is in the same spirit of {that} of \cite[Lemma 2.1]{Ng-CALR-frequency}. 
The existence of $v_\delta$ follows from the uniqueness of $v_{\delta}$ by using the limiting absorption principle, see,  e.g., \cite{Ng-WP}.  We now establish the uniqueness of $v_{\delta}$ by showing that $v_\delta = 0$ if $g_\delta =0$. Multiplying the equation of $v_\delta$ by $\bar v_\delta$ (the conjugate of $v_\delta$) and integrating by parts, we obtain 
\begin{equation*}
- \int_{\Omega} \langle A_\delta \nabla v_\delta, \nabla v_\delta \rangle  + k^2 \int_{\Omega} \Sigma_\delta |v_\delta|^2 + \int_{\partial \Omega} i k |v_\delta|^2  =0. 
\end{equation*}

Considering the imaginary part, and using the definition~\eqref{def_coeff} 
of $A_\delta$ and $\Sigma_\delta$, we have 
\begin{equation}\label{part1-Lem1}
v_\delta = 0 \mbox{ in } D. 
\end{equation}
This implies $v_\delta \big|_{D} = A_\delta \nabla v_\delta \big|_{D}  \cdot \nu = 0$ 
on $\partial D$; 
which yields, by the transmission conditions on $\partial D$,  
$$
v_\delta \big|_{\Omega \setminus D} = 
A  \nabla v_\delta \big|_{\Omega \setminus D}  \cdot \nu = 0 \mbox{ on } \partial D. 
$$
It follows from the unique continuation (see e.g. \cite{Protter60}) that $v_\delta = 0$ 
also in $\Omega \setminus D$. The proof of uniqueness is complete. 

We next establish \eqref{lem1-part1}  by  contradiction. Assume that 
there exists $(g_{\delta}) \subset [H^{1}(\Omega)]^*$ such that 
\begin{equation}\label{contradict-assumption}
\| v_\delta\|_{H^1(\Omega)} =1 \mbox{ and }\frac{1}{\delta} \Big| \int_{\Omega} g_\delta \bar v_\delta \Big|  + \|g_\delta \|_{[H^{1}(\Omega)]^*}^2 \to 0,
\end{equation}
as $\delta \to \hat  \delta \in [0, \delta_0]$. In fact, by contradiction these properties only hold for a sequence  $(\delta_n) \to \hat \delta$. However, for the simplicity of notation, we still use $\delta$ instead of $\delta_n$  to denote an element of such a sequence. 
We only consider the case   $\hat \delta = 0$; the case $\hat \delta > 0$ follows similarly. 
Without loss of generality, one may assume that $(v_\delta)$ converges to  $v_0$ strongly in $L^2(\Omega)$ and  weakly in $H^1(\Omega)$ for some $v_0 \in H^1(\Omega)$. Then, by \eqref{contradict-assumption},  
\begin{equation}\label{eq-v0-11}
\left\{\begin{array}{cl}
\dive(A_0 \nabla v_0) + k^2 \Sigma_0 v_0 =0 &  \mbox{ in } \Omega, \\[6pt]
A \nabla v_0 \cdot \nu - i k v_0  = 0 & \mbox{ on } \partial \Omega. 
\end{array}\right.
\end{equation}
Multiplying the equation of $v_\delta$ by $\bar v_\delta$ 
and integrating by parts, we obtain 
\begin{equation}\label{lem1-to0}
- \int_{\Omega} \langle A_\delta \nabla v_\delta, \nabla v_\delta \rangle  + k^2 \int_{\Omega} \Sigma_\delta |v_\delta|^2 + \int_{\partial \Omega} i k |v_\delta|^2  = \int_{\Omega} g_\delta \bar v_\delta. 
\end{equation}
Considering the imaginary part of \eqref{lem1-to0} and using \eqref{contradict-assumption}, we have 
\begin{equation}\label{lem1-to1}
\lim_{\delta \to 0} \Big( \| \nabla v_\delta\|_{L^2(D)} 
+ \| v_\delta \|_{L^2(D)} + \| v_\delta \|_{L^2({\partial \Omega})} \Big) = 0. 
\end{equation}
This implies $v_0 = 0$ in $D$
{and that $v_0 = 0$ on $\partial \Omega$}. 
As in the proof of uniqueness, we derive that $v_0 = 0$ in $\Omega$. Since $v_\delta \to v_0$ in $L^2(\Omega)$, it follows that 
\begin{equation}\label{lem1-to2}
\lim_{\delta \to 0} \|v_\delta \|_{L^2(\Omega)} = 0. 
\end{equation}
Considering the real part of \eqref{lem1-to0} and using \eqref{contradict-assumption}, \eqref{lem1-to1}, and \eqref{lem1-to2}, we obtain 
\begin{equation}\label{lem1-to2}
\lim_{\delta \to 0} \|\nabla v_\delta \|_{L^2(\Omega \setminus D)} = 0. 
\end{equation}
Combining \eqref{lem1-to0}, \eqref{lem1-to1}, and \eqref{lem1-to2} yields 
\begin{equation*}
\lim_{\delta \to 0} \|v_\delta \|_{H^1(\Omega)} = 0:  
\end{equation*}
which contradicts \eqref{contradict-assumption}. The proof is complete.  \proofend

\begin{remark} \label{rem-D-0} \fontfamily{m} \selectfont   In the case $k= 0$, the result in Lemma~\ref{lem1} also holds for zero Dirichlet boundary condition in which ${g}$
{may only be required to be in $L^2(\Omega)$}. The proof {follows the same lines}. 
\end{remark}

The following standard result  
is {repeatedly used} in this paper: 
\begin{lemma} \label{lem2} 
Let $d  =2 , 3$, $k \ge  0$. Let $D, \, V,  \,\Omega$ be smooth bounded connected open subsets 
of $\mR^d$ such that $D \subset \subset \Omega$, $\partial D \subset V \subset \Omega$. 
Let $A$ be a  matrix-valued function  and 
$\Sigma$ be a complex function{,} 
both defined in $\Omega${,} such that
\begin{equation*}
A \textrm{ is {\bf uniformly elliptic} in $\Omega$ and } \Sigma \in L^\infty(\Omega) 
\mbox{ with } \Im(\Sigma) \ge 0 \mbox{ and } \Re(\Sigma) \ge c  >0,  
\end{equation*}
for some constant $c$. Assume that $A \in C^1(\Omega  \setminus D)$ and $A \in C^1(V \cap \bar D)$.  Let $g  \in L^2(\Omega)$ and in the case $k = 0$ assume in addition that $\int_{\Omega} g = 0$. 
There exists a unique solution $v  \in H^1(\Omega)$ if $k > 0$ 
(respectively $v \in H^1_{m}(\Omega)$ if $k = 0$) to the system 
\begin{equation*}
\left\{\begin{array}{cl}
\dive(A \nabla v) + k^2 \Sigma v = g &  \mbox{ in } \Omega, \\[6pt]
A \nabla v \cdot \nu - i k v   = 0 & \mbox{ on } \partial \Omega. 
\end{array}\right.
\end{equation*}
Moreover, 
\begin{equation}\label{lem2-stability}
\| v\|_{H^1(\Omega)} \le C \| g\|_{L^2(\Omega)} \quad \mbox{ and } \quad \| v\|_{H^2({V} \setminus D)} \le C \| g\|_{L^2(\Omega)},  
\end{equation}
for some positive constant $C$ independent of $f$. 
\end{lemma} 

\noindent{\bf Proof.} The existence, uniqueness, and the first inequality of \eqref{lem2-stability} 
{follow}
from the Fredholm theory by the uniform ellipticity of $A$ in $\Omega$ and the boundary condition used. The second inequality of \eqref{lem2-stability} can be obtained by 
{Nirenberg's method of difference quotients} 
(see, e.g., \cite{BrAnalyse1}) using the smoothness assumption of $A$ 
and the boundedness of $\Sigma$. The details are left to the reader. \proofend 

\subsection{A toy problem} \label{sect-toy1} 
In this section, we consider a toy problem for tuned superlensing using HMMs,
in which the geometry is rectangular.  Given three positive constants $l$, $L$ and $T$,
we define \footnote{Letters $c$, $l$, $r$ stand for center, left, and right.}
\begin{equation*}
\R = [-l, L] \times [0, 2 \pi], \quad \R_{l}  = [-l, 0] \times  [0, 2 \pi], 
\quad  \R_c = [0, T] \times [0, 2 \pi], \quad \R_r = [T, L] \times [0, 2 \pi]. 
\end{equation*}
Denote 
\begin{equation*}
\Gamma: = \partial \R, \quad \Gamma_{c, 0} = \{ 0 \} \times [0, 2 \pi], \quad \mbox{ and } 
\quad \Gamma_{c, T} = \{ T \} \times [0, 2 \pi]. 
\end{equation*}
Let $a$ be a uniformly elliptic matrix-valued function defined in $\R_l \cup \R_r$. 
We set 
\begin{equation*}
a_{\delta} = \left(
\begin{array}{cc}
1 - i \delta   &  0     \\[6pt]
 0 &     -1 -  i \delta
\end{array}
\right), 
\end{equation*}
and define 
\begin{equation*}
A_{\delta} = \left\{\begin{array}{cl} a & \mbox{ in } \R_l \cup \R_r, \\[6pt]
a_{ \delta}  & \mbox{ in } \R_c,  
\end{array} \right. 
\end{equation*}
so that the superlensing device  occupies the region $\R_c$.
For $f \in L^2(\R)$ with $\supp f \cap \R_c = \O$, let $u_\delta \in H^1_0(\R)$ be the unique solution to the equation 
\begin{equation}\label{u-delta-E}
\dive(A_\delta \nabla u_\delta) = f \mbox{ in } \R. 
\end{equation}
Assume that $\| u_\delta \|_{H^1(\R)}$ is bounded as $\delta \to 0$.
Then, up to a subsequence, $u_\delta$ converges weakly to some $u_0 \in H^1_0(\R)$. 
It is clear that $u_0$ is a solution to 
\begin{equation}\label{eq-key}
\dive(A_0 \nabla u_0) = f \mbox{ in } \R.  
\end{equation}
More precisely, $u_0 \in H^1_0(\R)$ satisfies \eqref{eq-key} if and only if 
$u_0$ satisfies the elliptic-hyperbolic system
\begin{equation*}
\dive(a \nabla u_0) = f  \mbox{ in } \R_l \cup \R_r   \quad \mbox{ and } \quad \partial_{x_1 x_1}^2 u_0 - \partial_{x_2 x_2}^2 u_0 = f \mbox{ in } \R_c, 
\end{equation*}
and the transmission conditions
\begin{eqnarray*}
\left\{ \begin{array}{lcl}
u_0 \big|_{\R_l} &=& u_0 \big|_{\R_c}
\\[6pt]
\partial_{x_1} u_0 \big|_{\R_l} &=& \partial_{x_1} u_0 \big|_{\R_c},
\end{array}
\right.
 \textrm{on}\; \Gamma_{c, 0}  \quad 
\mbox{ and } \quad 
\left\{ \begin{array}{lcl}
u_0 \big|_{\R_r} &=& u_0 \big|_{\R_c}
\\[6pt]
\partial_{x_1} u_0 \big|_{\R_r} &=& \partial_{x_1} u_0 \big|_{\R_c},
\end{array} \right.
\mbox{ on }\; \Gamma_{c, T}.
\end{eqnarray*}
This problem is ill-posed: in general, there is no solution in $H^1_0(\R)$,
and so, $\| u_\delta \|_{H^1(\R)} \to + \infty$, as $\delta \to 0$. 
Nevertheless,  for some special choices of $T$, discussed below, 
the problem is well-posed and its solutions have peculiar properties.  
\medskip 

To describe them, we introduce an ``effective domain" $\R_T = [-l, L - T] \times [0, 2 \pi]$ and 
\begin{equation*}
\hat A (x_1, x_2), \hat f(x_1, x_2) = \left\{\begin{array}{cl} a(x_1, x_2), f(x_1, x_2) & \mbox{ in } \R_l \\[6pt]
a(x_1+T, x_2), f(x_1+T, x_2) & \mbox{ in } \R_T \setminus \R_l. 
\end{array}\right.
\end{equation*}
In what follows, we assume that $\hat A \in C^2(\overline{\R_T})$.

\begin{proposition} \label{pro1} 
Let $0 < \delta < 1$, $f \in L^2 (\R)$, and $u_\delta \in H^1_0(\R)$ be the unique solution of \eqref{u-delta-E}. 
Assume that $T \in  2 \pi \mN_+$ and $\sup f \cap \R_c = \O$.  Then
\begin{equation}\label{pro1-1-part1}
\| u_\delta \|_{H^1} \le C \| f\|_{L^2(\R)} \quad \mbox{ and } \quad  
u_\delta  \to u_0 \quad\textrm{strongly in}\; H^1(\R), 
\end{equation}
where  $u_0 \in H^1_0 (\R)$ is  the unique solution of \eqref{u-delta-E} with $\delta = 0$
and $C$ is a positive constant independent of $\delta$ and $f$. We also have 
\begin{equation*}
u_0(x_1, x_2) = \left\{ \begin{array}{lc} \hat u(x_1, x_2) &  \mbox{ in } \R_l, \\[6pt]
\hat u(x_1 - T, x_2) & \mbox{ in } \R_r,  
\end{array}\right.
\end{equation*}
where $\hat u \in H^1_0(\R_T)$ is the unique solution to the equation 
\begin{equation}\label{def-U-toy}
\dive (\hat A \nabla \hat u ) = \hat f \mbox{ in } \R_T,
\end{equation}
\end{proposition}

\begin{remark}  \fontfamily{m} \selectfont  
It follows from Proposition~\ref{pro1} that $u_0$ can be computed as if the structure in $\R_c$ had disappeared. 
This phenomenon is similar to that in the Veselago setting: superlensing occurs. 
\end{remark}

\noindent{\bf Proof.} The proof of Proposition~\ref{pro1} 
is in the spirit of the approach used by the second author 
in ~\cite{Ng-Complementary} to deal with negative index materials. 
The key point is to construct the unique solution $u_0$ to the limiting problem 
appropriately and then obtain estimates on $u_\delta$ by studying the difference 
$u_\delta - u_0$.

We first  construct a solution $u_0 \in H^1_0(\R)$ to \eqref{u-delta-E} with $\delta =0$. Since $\hat A \in C^2(\overline{\R_T})$ and since $f \in L^2(\R)$, 
the regularity theory for elliptic equations (see,  e.g., \cite[3.2.1.2]{Grisvard}) implies that $\hat u \in H^2(\R)$
and 
\begin{equation}\label{pro1-est1}
\|\hat u \|_{H^2(\R)}\le C \|f \|_{L^2(\R)}. 
\end{equation}
Here and in what follows in this proof, $C$ denotes a positive constant independent of $f$ and $\delta$. It follows that $\hat{u}(0, x_2) \in H^1(\Gamma_{c, 0})$ 
and $\partial_1 \hat{u}(0, x_2) \in L^2(\Gamma_{c, 0})$. Interpretting $x_1$ and $x_2$ as respectively time and space variables in the rectangle $\R_c$,
 we  seek a solution $v \in C\big([0, T]; H^1_0(0, 2 \pi) \big) \cap C^1([0, T]; L^2(0, 2 \pi))$ 
of  the wave equation 
\begin{equation}\label{def-uc-toy}
\partial^2_{x_1 x_1} v - \partial^2_{x_2 x_2} v = 0 \quad \mbox{ in } \R_c,
\end{equation}
with zero boundary condition, i.e., $v = 0 $ on $\Gamma \cap \partial \Omega_c$, 
and  the following initial  conditions 
\begin{equation*}
v (0, x_2) = \hat u (0, x_2) \quad \mbox{ and } \quad \partial_{x_1} v (0, x_2) = \partial_{x_1} \hat u \big|_{\R_l} (0, x_2). 
\end{equation*}
Existence and uniqueness of $v$ follow from the standard theory of the wave equation by taking into account the regularity information in \eqref{pro1-est1}. We also have, for $0\le x_1 \le T$,  
\begin{align}\label{pro1-est2}
\int_{0}^{2 \pi} |\partial_{x_1} v(x_1, x_2)|^2 + |\partial_{x_2} v(x_1, x_2)|^2 \, d x_2 = &  \int_{0}^{2 \pi} |\partial_{x_1} v(0, x_2)|^2 + |\partial_{x_2} v(0, x_2)|^2 \, d x_2 \nonumber \\[6pt]
= &  \int_{0}^{2 \pi} |\partial_{x_1} \hat u \big|_{\R_l}(0, x_2)|^2 + |\partial_{x_2} \hat u (0, x_2)|^2 \, d x_2.
\end{align}
Furthermore, one  can represent $v$ in $\R_c$ in the form 
\begin{equation}\label{form-uc-toy}
v(x_1, x_2) = \sum_{n =  1}^\infty  \sin (n x_2) \big[a_n \cos(n x_1) + b_n \sin( n x_1)\big],
\end{equation}
where $a_n, b_n \in \mR$ are determined by the initial conditions satisfied by $v$ at $x_1= 0$. 
Since $T \in 2 \pi \mN$, it follows that 
\begin{equation}\label{key-ob-toy}
v(0, \cdot ) = v(T, \cdot ) \quad \mbox{ and } \quad \partial_{x_1} v(0, \cdot ) = \partial_{x_1} v(T, \cdot ) \mbox{ in }  [0, 2 \pi],
\end{equation} 
for any initial conditions, and hence for any $f$ with $\supp f \cap \R_c = \O$. 
Define 
\begin{equation}\label{def-u-toy}
u_0(x_1, x_2) =  \left\{ \begin{array}{lc} \hat u(x_1, x_2) &  \mbox{ in } \R_l, \\[6pt]
v(x_1, x_2) & \mbox{ in } \R_c, \\[6pt]
\hat u (x_1 - T, x_2) & \mbox{ in } \R_r. 
\end{array}\right.
\end{equation}
It follows from \eqref{def-U-toy}, \eqref{def-uc-toy},  and \eqref{key-ob-toy} 
that $u_0 \in H^1_0(\Omega)$ is a solution to \eqref{u-delta-E} with $\delta = 0$; moreover, by \eqref{pro1-est1} and \eqref{pro1-est2}, 
\begin{equation}\label{pro1-est3}
\| u_0\|_{H^1(\R)} \le C \|f \|_{L^2(\R)}. 
\end{equation}

We next establish the uniqueness of $u_0$. 
Let $w_0 \in H^1_0(\Omega)$ be a solution to \eqref{u-delta-E} with $\delta = 0$. 
Since $w_0$ can be represented as in \eqref{form-uc-toy} in $\R_c$, we obtain 
\begin{equation*}
w_0(0, \cdot ) = w_0(T, \cdot ) \quad  \mbox{ and } 
\quad \partial_{x_1} w_0(0, \cdot ) = \partial_{x_1} w_0(T, \cdot ) \mbox{ in }  [0, 2 \pi]. 
\end{equation*}
We can thus define for $(x_1, x_2)$ in $\R_T$  
\begin{equation*}
\hat w (x_1, x_2) = \left\{\begin{array}{cl} w_{0}(x_1, x_2) & \mbox{ in } \R_l, \\[6pt]
w_{0}(x_1 - T, x_2) & \mbox{ otherwise},
\end{array}\right.
\end{equation*}
which is a solution to~\eqref{def-U-toy}. By uniqueness, it follows that
$\hat{w} \equiv \hat{u}$ in $\R_T$,
and \eqref{def-u-toy} shows that $w_0 \equiv u_0$ in $\R$.
\medskip 

Finally, we establish \eqref{pro1-1-part1}. Define 
\begin{equation}\label{def-vd}
v_\delta = u_\delta - u_0 \;\mbox{ in } \R. 
\end{equation} 
We have 
\begin{align*}
\dive(A_\delta \nabla v_\delta) = &  \dive(A_\delta \nabla u_\delta) - \dive(A_\delta \nabla u_0) \\[6pt]
=  & \dive(A_\delta \nabla u_\delta) - \dive(A_0 \nabla u_0) +  \dive(A_0 \nabla u_0) - \dive(A_\delta \nabla u_0)
\mbox{ in } \R. 
\end{align*}
It follows that  $v_\delta \in H^1_0(\Omega)$ is the solution to
\begin{equation}\label{hahaTT}
\dive(A_\delta \nabla v_\delta) =  \dive (i \delta  \mathds{1}_{\R_c} \nabla u_0) \mbox{ in } \R. 
\end{equation}
As in \eqref{lem1-part2} in Lemma~\ref{lem1}, we obtain from \eqref{pro1-est3} that 
\begin{equation}\label{pro1-est4}
\| v_\delta \|_{H^1(\R)} \le \frac{C}{\delta} \|  \delta  \nabla u_0  \|_{L^2(\R_c)} \le  C \| f\|_{L^2(\R)};  
\end{equation}
which implies the first inequality of  \eqref{pro1-1-part1}. 
As in \eqref{lem1-part1} in Lemma~\ref{lem1}, we deduce from \eqref{pro1-est3}, \eqref{hahaTT}, and \eqref{pro1-est4} that 
\begin{equation} \label{pro1-est5}
\|u_\delta - u_0\|_{H^1(\R)}^2 = \|v_\delta\|_{H^1(\R)}^2 
\le C \left| \int_{\R_c} i \nabla u_0 \nabla v_\delta \right| \to 0, 
\end{equation}
as $v_\delta$ converges weakly to $0$ in $H^1(\R)$.  The proof is complete. \proofend


\subsection{Tuned superlensing using HMMs}

In this section we consider a superlens of the form~\eqref{scheme1},
with the constraint~\eqref{cond-r1-r2}.
We establish a more general version of Theorem~\ref{thm-main}
and its variants in the finite frequency regime. 
We then present another scheme in the same spirit, in which the superlens 
is strictly hyperbolic and not merely degenerately hyperbolic.  
\medskip

We first deal with a  situation in two dimensions. We consider a cylindrical lens, defined in $B_{r_2} \setminus B_{r_1}$ by 
a pair $(A^H, \Sigma^H)$ of the form
\begin{equation}\label{scheme1-k}
(A^H, \Sigma^H) = (\frac{1}{r}e_r \otimes e_r - r e_\theta \otimes e_\theta, 0) \mbox{ in } B_{r_2} \setminus B_{r_1}. 
\end{equation}
Assume that the region $B_{r_1}$ to be magnified
is characterized by a pair $(a, \sigma)$  of a matrix-valued function $a$ and a complex function $\sigma$ such that $a$ satisfies the standard condition mentioned in the introduction ($a$ is uniformly elliptic in $B_{r_1}$ and \eqref{smoothness} holds) and $\sigma$ satisfies the following standard conditions
\begin{equation}\label{cond-sigma}
 \sigma \in L^\infty(B_{r_1}), 
\;\textrm{with}\; \Im(\sigma) \ge 0 \mbox{ and } \Re(\sigma) \ge c  > 0, 
\end{equation}
for some constant  $c$. 

Taking loss into account, the overall medium is characterized by
\begin{equation}\label{def-AdS-thm2}
A_\delta, \Sigma_\delta  = \left\{\begin{array}{cl} I, 1 & \mbox{ in } \Omega \setminus B_{r_2}, \\[6pt]
A^H - i \delta I, \Sigma^H + i \delta & \mbox{ in } B_{r_2} \setminus B_{r_1}, \\[6pt]
a, \sigma & \mbox{ in } B_{r_1},
\end{array}\right.
\end{equation}
Given a (source) function $f \in L^2(\Omega)$ and given a frequency $k>0$, 
standard arguments show that there is a unique solution $u_\delta \in H^1(\Omega)$
to the system
\begin{equation}\label{sys-k-scheme1}
\left\{\begin{array}{cl}
\dive (A_\delta \nabla u_\delta) + k^2 \Sigma_\delta u_\delta = f & \mbox{ in } \Omega, \\[6pt]
\partial_{\nu} u_\delta - i k u_\delta = 0 &  \mbox{ on } \partial \Omega. 
\end{array}\right. 
\end{equation}
The following theorem describes the superlensing property of the superlensing device defined 
by~\eqref{scheme1-k}.

\begin{theorem}\label{thm-scheme1-k}
Assume $d=2$ and $k>0$. Let $0< \delta < 1$,  $\Omega$ be a smooth bounded connected open 
subset of $\mR^2$, and let $f \in L^2(\Omega)$. 
Let $0 < r_1 < r_2$ be such that \eqref{cond-r1-r2} holds, 
and assume that $B_{r_2} \subset \subset \Omega$, and $\supp f \subset \Omega \setminus B_{r_2}$.
Let $u_\delta \in H^1(\Omega)$ be the unique solution to \eqref{sys-k-scheme1}. 
We have 
\begin{equation}\label{thm-scheme1-k-statement1}
\|u_\delta  \|_{H^1(\Omega)} \le C  \| f\|_{L^2(\Omega)} \quad \mbox{ and } \quad 
u_\delta \rightarrow u_0 \quad\textrm{strongly in}\; H^1(\Omega), 
\end{equation}
where $u_0 \in H^1(\Omega)$ is 
the unique solution to~\eqref{sys-k-scheme1} with $\delta =0$ and $C$ is a positive constant independent of $f$ and $\delta$. 
Moreover, $u_0 = \hat u$ in $\Omega \setminus B_{r_2}$, 
where $\hat u$ is the unique solution to the system 
\begin{equation}\label{def-hA-hu-k}
\left\{\begin{array}{cl}
\dive(\hat A  \nabla \hat u) + k^2 \hat \Sigma \hat u = f &  \mbox{ in } \Omega\\[6pt]
\partial_\nu \hat u - i k \hat u  = 0 &  \mbox{ on } \partial \Omega, 
\end{array}\right. 
\end{equation}
where
\begin{equation*}
\hat A(x), \hat \Sigma (x)  = \left\{\begin{array}{cl}  I, 1 & \mbox{ in } \Omega \setminus B_{r_2},  \\[6pt]
\dsp a\Big( \frac{r_1 }{ r_2 }  x\Big),  \frac{r_1^2}{ r_2^2} \sigma \Big( \frac{r_1 }{ r_2 }  x\Big)& \mbox{ in } B_{r_2}.
\end{array}\right.
\end{equation*}
\end{theorem}

\noindent 
Since $f$ is arbitrary with support in $\Omega \setminus B_{r_2}$, 
it follows from the definition of $(\hat A, \hat \Sigma)$ that the object 
in $B_{r_1}$ is magnified by a factor $r_2/ r_1$. 
\medskip 


\noindent
{\bf Proof of Theorem~\ref{thm-scheme1-k}.} 
The proof is in the spirit of that of Proposition~\ref{pro1}: the main idea is to 
construct $u_0$ and then estimate $u_\delta - u_0$. 

We have 
\begin{equation}\label{thm-scheme1-k-part0}
\| \hat u\|_{H^1(\Omega)} \le C \| f\|_{L^2(\Omega)}. 
\end{equation}
Using \eqref{smoothness} and applying Lemma~\ref{lem2},  we derive  that $u \in H^2(\Omega \setminus B_{r_2})$ and 
\begin{equation}\label{thm-scheme1-k-part1}
\|\hat u \|_{H^2(\Omega \setminus B_{r_2})} \le C \| f\|_{L^2(\Omega)}. 
\end{equation}
Define a function $v$ in $B_{r_2} \setminus B_{r_1}$ by 
\begin{equation}\label{WE-thm1}
\partial^2_{rr} v  - \partial^2_{\theta \theta} v = 0, 
\quad v \;\textrm{is periodic with respect to $\theta$}, 
\end{equation}
and 
\begin{equation}\label{WE-thm1-bdry}
v(r_2, \theta) = \hat u(r_2, \theta)  \quad \mbox{ and } \quad 
\partial_r v (r_2, \theta) = \partial_r \hat u (r_2,  \theta) \big|_{\Omega \setminus B_{r_2}}  
\mbox{ for } \theta \in [0, 2 \pi]. 
\end{equation}
By considering \eqref{WE-thm1} as a Cauchy problem for the wave equation with periodic boundary conditions, 
in which  $r$ and $\theta$ are seen as a time and a space variable respectively, 
the standard theory shows that there exists a unique such 
$v(r, \theta) \in C\big([r_1, r_2]; H^1_{per}(0, 2 \pi) \big) \cap C^1([r_1, r_2]; L^2(0, 2 \pi))$. 
We also have, for $r_1 \le r \le r_2$,  
\begin{align}\label{thm-scheme1-k-est2}
\int_{0}^{2 \pi} |\partial_{r} v(r, \theta)|^2 + |\partial_{\theta} v(r, \theta)|^2 \, d \theta = &  \int_{0}^{2 \pi} |\partial_{r} v(r_2, \theta)|^2 + |\partial_{\theta} v(r_2, \theta)|^2 \, d \theta \nonumber \\[6pt]
= &  \int_{0}^{2 \pi} |\partial_{r} \hat u \big|_{\R_l}(r_2, \theta)|^2 + |\partial_{\theta} \hat u (r_2, \theta)|^2 \, d \theta; 
\end{align}
which yields,  by \eqref{thm-scheme1-k-part1},  
\begin{equation}\label{thm-scheme1-k-part2}
\| v\|_{H^1(B_{r_2} \setminus B_{r_1})} \le C \|f \|_{L^2(\Omega)}. 
\end{equation}
Moreover, $v$ can be represented in the form   
\begin{equation}\label{form-uc}
v(r, \theta) = \sum_{-\infty}^\infty \sum_{\pm} a_{n, \pm} e^{i(n r \pm n \theta)}  
\;\textrm{in}\; B_{r_2} \setminus B_{r_1}, 
\end{equation}
where $a_{n, \pm} \in \mC$. 
Since $r_2 - r_1 \in 2 \pi \mN_+$, it follows that
\begin{equation}\label{W-thm1}
v(r_1, \theta) = v(r_2, \theta) \quad \mbox{ and } \quad \partial_r v(r_1, \theta) = \partial_r v(r_2, \theta) \mbox{ for } \theta \in [0, 2\pi]. 
\end{equation}
Set  
\begin{equation}\label{def-u0-scheme1}
u_0 = \left\{ \begin{array}{cl} 
\hat u &  \mbox{ in } \Omega \setminus B_{r_2}, \\[6pt]
v & \mbox{ in } B_{r_2} \setminus B_{r_1}, \\[6pt]
\dsp \hat u \big(r_2 x/ r_1 \big) & \mbox{ in } B_{r_1}.  
\end{array}\right.
\end{equation}
It follows from \eqref{thm-scheme1-k-part0} and \eqref{thm-scheme1-k-part2} that 
\begin{equation}\label{thm-scheme1-k-part3}
\|u_0 \|_{H^1(\Omega)} \le C \|f \|_{L^2(\Omega)}. 
\end{equation}
We also have 
 \begin{equation}\label{part1-u0}
\dive (A_0 \nabla u_0) + k^2 \Sigma_0 u_0 = f 
\;\textrm{in}\; \Omega \setminus (\partial B_{r_1} \cup \partial B_{r_2}). 
\end{equation}
On the other hand, from \eqref{WE-thm1-bdry} and the definition of $\hat A$, we have
\begin{equation}\label{part2-u0}
[u_0] = [A_0 \nabla u_0 \cdot e_r ] = 0 \mbox{ on } \partial B_{r_2}
\end{equation}
and from \eqref{W-thm1},  we obtain 
\begin{equation}\label{part3-u0}
[u_0] = [ A_0 \nabla u_0 \cdot e_r ] = 0 \mbox{ on } \partial B_{r_1}. 
\end{equation}
A combination of \eqref{part1-u0}, \eqref{part2-u0}, and \eqref{part3-u0} yields that 
\begin{equation*}
\dive (A_0 \nabla u_0) + k^2 \Sigma_0 u_0 = f \mbox{ in } \Omega;  
\end{equation*}
which implies that  $u_0$ is a solution to \eqref{sys-k-scheme1} with $\delta = 0$. 
\medskip

We next establish the uniqueness of $u_0$. Let $w_0 \in H^1(\Omega)$ be a solution 
to \eqref{sys-k-scheme1} with $\delta = 0$. 
Since $w_0$ can be represented as in \eqref{form-uc} in $B_{r_2} \setminus B_{r_1}$, we have
\begin{equation}\label{w0-pro}
w_0(r_1, \theta) = w_0(r_2, \theta) \quad \mbox{ and } 
\quad \partial_r w_0(r_1, \theta) = \partial_r w_0(r_2, \theta) \mbox{ for } \theta \in [0, 2\pi]. 
\end{equation}
Define 
\begin{equation*}
\hat w(x)  = \left\{\begin{array}{cl} w_{0}(x) & \mbox{ in } \Omega \setminus B_{r_2}, \\[6pt]
\dsp w_{0}\big(r_1 x/ r_2 \big) & \mbox{  in } B_{r_2}. 
\end{array}\right.
\end{equation*}
It follows from \eqref{w0-pro} that $\hat w \in H^1(\Omega)$. It is easy to verify that $\hat w$ is a solution of \eqref{def-hA-hu}. Hence $\hat w = \hat u$; which yields $w_0 = u_0$. 
\medskip

We finally establish \eqref{thm-scheme1-k-statement1}. Set 
\begin{equation}\label{def-vd2}
v_\delta = u_\delta - u_0 \mbox{ in } \Omega.  
\end{equation} 
Then $v_\delta \in H^1(\Omega)$ and satisfies 
\begin{equation*}
\dive(A_\delta \nabla v_\delta) + k^2 \Sigma_\delta v_\delta =  \dive (i \delta   \mathds{1}_{B_{r_2} \setminus B_{r_1}} \nabla u_0) - i \delta k^2 \mathds{1}_{B_{r_2} \setminus B_{r_1}} u_0 \mbox{ in } \Omega. 
\end{equation*}
and 
$$
\partial_\nu v_\delta - i k v_\delta = 0 \mbox{ on } \partial \Omega. 
$$
Applying \eqref{lem1-part2} of Lemma~\ref{lem1}, we obtain from \eqref{thm-scheme1-k-part3} that 
$$
 \| v_\delta \|_{H^1(\Omega)} \le C \| \nabla u_0\|_{L^2(\Omega)},
$$ 
which is the inequality in \eqref{thm-scheme1-k-statement1}.  
Applying \eqref{lem1-part1} of Lemma~\ref{lem1}, we derive from \eqref{thm-scheme1-k-part3} that 
$$
\|u_\delta - u_0\|_{H^1(\Omega)}^2 = \|v_\delta\|_{H^1(\Omega)}^2 
\le C \left| \int_{B_{r_2} \setminus B_{r_1}} i \nabla u_0 \nabla v_\delta \right| 
\;\to\; 0, 
$$ 
which completes the proof. 
\proofend


\begin{remark}  \fontfamily{m} \selectfont  
The  proof of Theorem~\ref{thm-main} when  $A^H$ is given by \eqref{scheme1}-\eqref{cond-r1-r2} in two dimensions is similar to the one of Theorem~\ref{thm-scheme1-k}. The details are left to the reader.  
\end{remark}

In the rest of this section, we consider another construction, 
for the three dimensional finite frequency case, in which the superlens 
is made of (strictly) hyperbolic metamaterials. 
Instead of \eqref{scheme1-k}, the superlens is now defined by
\begin{equation}\label{scheme1-k-2}
(A^H, \Sigma^H) = \Big(\frac{1}{r^2}e_r \otimes e_r -  (e_\theta \otimes e_\theta + e_\varphi \otimes e_\varphi) ,   \frac{1}{4 k^2r^2} \Big) \mbox{ in } B_{r_2} \setminus B_{r_1}. 
\end{equation}
Note that $\Sigma^H$ now depends on $k$.  

\begin{theorem}\label{thm-scheme1-k-2} 
Let $d=3$, $k>0$, $0 < \delta < 1$,  and let $\Omega$ be a smooth bounded connected open subset of $\mR^3$. 
Let $f \in L^2(\Omega)$.
Fix $0 < r_1 < r_2$, such that $r_2 - r_1 \in  4 \pi \mN_+$,  $B_{r_2} \subset \subset \Omega$ and $\supp f \subset \Omega \setminus B_{r_2}$. 
Let $u_\delta \in H^1(\Omega)$ be the unique solution of \eqref{sys-k-scheme1}, 
where $(A^H, \Sigma^H) $ is given by \eqref{scheme1-k-2}.   
We have 
\begin{equation}\label{thm-scheme1-k-2-statement1}
\|u_\delta  \|_{H^1(\Omega)} \le C  \| f\|_{L^2(\Omega)} \quad \mbox{ and } \quad 
u_\delta \rightarrow u_0 \quad \textrm{strongly in}\; H^1(\Omega), 
\end{equation}
where $u_0 \in H^1_{m}(\Omega)$ is 
the unique solution to~\eqref{sys-k-scheme1} with $\delta =0$ and $C$ is a positive constant independent of $f$ and $\delta$. 
Moreover, $u_0 = \hat u$ in $\Omega \setminus B_{r_2}$ where $\hat u$ 
is the unique solution to the system \eqref{def-hA-hu}, where
\begin{equation*}
\hat A(x), \hat \Sigma (x)  = \left\{\begin{array}{cl}  I, 1 & \mbox{ in } \Omega \setminus B_{r_2},  \\[6pt]
\dsp \frac{r_1}{r_2}  a\Big( \frac{r_1 }{ r_2 }  x\Big),  \frac{r_1^3}{ r_2^3} \sigma \Big( \frac{r_1 }{ r_2 }  x\Big)& \mbox{ in } B_{r_2}.
\end{array}\right.
\end{equation*} 
\end{theorem}
\medskip

From the definition of $(A^H, \Sigma^H)$ in \eqref{scheme1-k-2}, one derives that  if $u$ is a solution to the equation $\dive(A^H \nabla u) + k^2 \Sigma^H u = 0$ in $B_{r_2} \setminus B_{r_1}$ then 
\begin{equation*}
\partial^2_{rr} u - \Delta_{\partial B_1} u + \frac{1}{4} u = 0 \mbox{ in } B_{r_2} \setminus B_{r_1}. 
\end{equation*}
This equation plays a similar role as the wave equation~\eqref{wave1-ole}. 
The proof of Theorem~\ref{thm-scheme1-k-2} below shows that 
\begin{eqnarray*}
u(r_1 \hat{x}) = u(r_2 \hat{x}) 
&\quad \textrm{and} \quad&
\partial_r u(r_1 \hat{x}) = \partial_r u(r_2 \hat{x}),
\quad \textrm{for}\;  \hat{x} \in \partial B_1. 
\end{eqnarray*}
The same strategy as that used for proving Theorem~\ref{thm-scheme1-k}, then 
leads to the above conclusion.
\medskip 

\noindent{\bf Proof.} We have 
$$
\|\hat u\|_{H^1(\Omega)} \le C \| f\|_{L^2(\Omega)}
$$ 
and,  by \eqref{smoothness} and Lemma~\ref{lem2}, 
$$
\|\hat u\|_{H^2(\Omega \setminus B_{r_2})} \le C \| f\|_{L^2(\Omega)}. 
$$ 
For $n \ge 0$ and $-n \le m \le n$, let $Y^m_n$ denote the spherical harmonic function 
of degree $n$ and of order $m$, which satisfies
\begin{eqnarray*}
\Delta_{\partial B_1} Y^m_n + n(n+1) Y^m_n 
&=& 0 \quad\textrm{on}\; \partial B_1. 
\end{eqnarray*}
Since the family $\big(Y^m_n \big)$ is dense in $L^2(\partial B_1)$, 
any solution $v \in H^1(B_{r_2} \setminus B_{r_1})$ to the equation 
\begin{equation}\label{haha}
\partial^2_{rr} v- \Delta_{\partial B_1} v +  \frac{1}{4} v = 0 \mbox{ in } B_{r_2} \setminus B_{r_1}, 
\end{equation}
can be represented in the form
\begin{eqnarray} \label{form_v}
v(x) &=& 
\sum_{n = 1}^\infty \sum_{m =-n}^n \sum_{\pm} 
a_{nm, \pm} e^{\pm i \lambda_n r} Y_m^n(\hat{x}),
\quad\quad x \in B_{r_2} \setminus B_{r_1}, 
\end{eqnarray}
where $\lambda_n = (n+1/2)$,  $r= |x|$ and $\hat{x} = \frac{x}{|x|}$.
Note that the 0-order term in \eqref{haha} has been chosen in $B_{r_2} \setminus B_{r_1}$
so that the dispersion relation writes
\begin{equation*}
\lambda_n^2 = n(n+1) + \displaystyle\frac{1}{4} = (n + \displaystyle\frac{1}{2})^2,
\end{equation*}
which implies that all the terms $e^{\pm i\lambda_n r}$ in~\eqref{form_v}, and thus $v$, 
are $4 \pi$-periodic functions of $r$.
Since $r_2 - r_2  \in 4 \pi \mN_+$, it follows that
\begin{equation*}
v(r_1 \hat x) = v(r_2 \hat x) \quad \mbox{ and } 
\quad \partial_r v(r_1 \hat x) = \partial_r v(r_2 \hat x) 
\quad \textrm{for}\;  \hat x \in \partial B_1. 
\end{equation*}
The conclusion follows as in the proof of Theorem~\ref{thm-scheme1-k} by noting that $u_0$ is also given by \eqref{def-u0-scheme1}. The details are left to the reader. 
\proofend


\section{Superlenses using  HMMs via complementary property}\label{sect-complementary}

In this section, we consider a lens with coefficients $A^H$
of the form \eqref{scheme2} in the finite frequency regime,
and we prove a superlensing result. This proof can be easily adapted
to obtain the conclusion of  Theorem~\ref{thm-main}, which corresponds to the quasistatic case.
\medskip

The superlensing device characterized by $(A^H, \Sigma^H)$ defined in $B_{r_2} \setminus B_{r_1}$ in the finite frequency regime is given by:   
\begin{equation}\label{scheme2-k}
(A^H, \Sigma^H)= 
\left\{ \begin{array}{cl} \dsp \Big( \frac{1}{r^{d-1}}e_r \otimes e_r 
- r^{3-d} (I - e_r \otimes e_r), \frac{1}{r^2} \Big)  
& \mbox{ in } B_{r_2} \setminus B_{r_m}, \\[6pt]
\dsp \Big(-  \frac{1}{r^{d-1}} e_r \otimes e_r 
+ r^{3-d} (I - e_r \otimes e_r), - \frac{1}{r^2} \Big)  
& \mbox{ in }  B_{r_m} \setminus B_{r_1}. 
\end{array} \right.
\end{equation}
Recall that  
\begin{equation*}
r_{m} = (r_1 + r_2 )/ 2. 
\end{equation*}
It will be clear below, that the choice $\Sigma^H = 1/r^2$ in $B_{r_2} \setminus B_{r_m}$ 
and $-1/r^2$ in  $B_{r_m} \setminus B_{r_1}$ is just a matter of simplifying the presentation. 
Any real-valued pair 
$(\tilde \sigma_1/r^2, \tilde \sigma_2/r^2) \in L^\infty(B_{r_m} \setminus B_{r_1}) 
\times L^\infty(B_{r_2} \setminus B_{r_m})$ 
which satisfies 
$$
\tilde \sigma_2(x) = - \tilde \sigma_1 \big([|x| - r_m]x / |x| \big) 
$$
is admissible. 
Assume that the region $B_{r_1}$ to-be-magnified contains a medium 
characterized by a pair  $(a, \sigma)$ of a matrix-valued function $a$ and a complex functions $\sigma$ such that 
$a$ satisfies the standard condition mentioned in the introduction ($a$ is uniformly elliptic in $B_{r_1}$ and \eqref{smoothness} holds) and $\sigma$ satisfies  \eqref{cond-sigma}.

 In the presence of the superlensing device and the object, the medium in $\Omega$ with the loss is characterized by
\begin{eqnarray}\label{def-AdS-thm2}
A_\delta, \Sigma_\delta  
&=&
\left\{\begin{array}{cl} I, 1 & \mbox{ in } \Omega \setminus B_{r_2}, \\[6pt]
A^H - i \delta I, \Sigma^H + i \delta & \mbox{ in } B_{r_2} \setminus B_{r_1}, \\[6pt]
a, \sigma & \mbox{ in } B_{r_1}.
\end{array}\right.
\end{eqnarray}
Given a source $f \in L^2(\Omega)$ and given a frequency $k>0$, 
the electromagnetic field $u_\delta$ is the unique solution to the system 
\begin{equation}\label{sys-k-scheme2}
\left\{\begin{array}{cl}
\dive (A_\delta \nabla u_\delta) + k^2 \Sigma_\delta u_\delta = f & \mbox{ in } \Omega, \\[6pt]
\partial_{\nu} u_\delta - i k u_\delta = 0 &  \mbox{ on } \partial \Omega. 
\end{array}\right.
\end{equation}
The superlensing property of the device~\eqref{scheme2-k} 
is given by the following theorem: 

\begin{theorem}\label{thm-scheme2-k} 
Let $d=2, 3$, $k>0$, $\Omega$ be a smooth bounded connected open subset of $\mR^d$, 
and let $f \in L^2(\Omega)$. 
Fix $0 < r_1 < r_2$, such that \eqref{cond-r1-r2} holds and assume 
that $B_{r_2} \subset \subset \Omega$ and  $\supp f \subset \Omega \setminus B_{r_2}$. 
Let $u_\delta \in H^1(\Omega)$ be the unique solution to~\eqref{sys-k-scheme2}. 
We have 
\begin{equation}\label{thm-scheme2-k-statement1}
\|u_\delta  \|_{H^1(\Omega)} \le C  \| f\|_{L^2(\Omega)} \quad \mbox{ and } \quad 
u_\delta \rightarrow u_0 \quad \textrm{strongly in}\; H^1(\Omega),
\end{equation}
where $u_0 \in H^1(\Omega)$ is 
the unique solution to~\eqref{sys-k-scheme2} with $\delta =0$ and $C$ is a positive constant independent of $f$ and $\delta$. 
Moreover, $u_0 = \hat u$ in $\Omega \setminus B_{r_2}$, 
where $\hat u_\delta$ is the unique solution to the system 
\begin{equation}\label{def-hA-hu-2}
\left\{\begin{array}{cl}
\dive(\hat A  \nabla \hat u) + k^2 \hat \Sigma \hat u = f &  \mbox{ in } \Omega\\[6pt]
\partial_\nu \hat u - i k \hat u  = 0 &  \mbox{ on } \partial \Omega,
\end{array}\right. 
\end{equation}
where
\begin{equation*}
\hat A(x), \hat \Sigma (x)  = \left\{\begin{array}{cl}  I, 1 & \mbox{ in } \Omega \setminus B_{r_2},  \\[6pt]
\dsp \frac{r_1^{d-2}}{r_2^{d-2}} a\Big( \frac{r_1 }{ r_2 }  x\Big),  \frac{r_1^{d}}{r_2^{d}} \sigma \Big( \frac{r_1 }{ r_2 }  x\Big)& \mbox{ in } B_{r_2}.
\end{array}\right.
\end{equation*}
\end{theorem}
\medskip

Since $f$ is arbitrary with support in $\Omega \setminus B_{r_2}$, 
it follows from the definition of $\hat A$ that the object in $B_{r_1}$ is 
magnified by a factor $r_2/ r_1$. 
We emphasize again that no condition is imposed on $r_2 - r_1$.  
\medskip

\noindent{\bf Proof.} 
The proof is in the spirit of that of Theorem~\ref{thm-scheme1-k}: 
the main idea is to construct $u_0$, solution to~\eqref{sys-k-scheme2} for 
$\delta = 0$, from $\hat u$ via reflection as discussed in the introduction, 
and then to estimate $u_\delta - u_0$. 

We have 
\begin{equation}\label{thm-scheme2-k-est0}
\| \hat u\|_{H^1(\Omega)} \le C \| f\|_{L^2(\Omega)}, 
\end{equation}
and,  by \eqref{smoothness} and Lemma~\ref{lem2}, 
\begin{equation}\label{thm-scheme2-k-est1}
\| \hat u\|_{H^2(\Omega \setminus B_{r_2})} \le C \| f\|_{L^2(\Omega)}. 
\end{equation}
Define  $v$  in $B_{r_2} \setminus B_{r_m}$ as follows 
\begin{equation}\label{WE-thm2}
\partial^2_{rr} v  - \Delta_{\partial B_1} v +k^2 v = 0 \mbox{ in } B_{r_2} \setminus B_{r_m}
\end{equation}
and 
\begin{equation}\label{WE-thm2-bdry}
v = \hat u  \quad \mbox{ and } \quad \partial_r v = \partial_r \hat u \big|_{\Omega \setminus B_{r_2}}  \mbox{ on } \partial B_{r_2}. 
\end{equation}
Consider \eqref{WE-thm2} and \eqref{WE-thm2-bdry} as a Cauchy problem for the wave equation defined on the manifold $\partial B_1$ for which $r$ plays as a time variable. 
By the standard theory for the wave equation, there exists a unique such  
$v \in C\big([r_m, r_2]; H^1(\partial B_1) \big) \cap C^1([r_m, r_2]; L^2(\partial B_1))$. We also have 
\begin{multline}\label{thm-scheme2-k-est2}
\int_{\partial B_1} |\partial_{r} v(r, \xi)|^2 + |\nabla_{\partial B_1} v(r, \xi)|^2 + k^2 |v(r, \xi)|^2 \, d \xi \\[6pt]
=   \int_{\partial B_1} |\partial_{r} v(r_2, \xi)|^2 + |\nabla_{\partial B_1} v(r_2, \xi)|^2 + k^2 |v(r_2, \xi)|^2 \, d \xi \\[6pt] =   \int_{\partial B_1} |\partial_{r} \hat u \big|_{\Omega \setminus B_{r_2}} (r_2, \xi)|^2 + |\nabla_{\partial B_1} \hat u (r_2, \xi)|^2 + k^2 |\hat u (r_2, \xi)|^2 \, d \xi. 
\end{multline}
It follows that $v \in H^1(B_{r_2} \setminus B_{r_m})$ and 
\begin{equation}\label{thm-scheme2-k-est3}
\| v\|_{H^1(B_{r_2} \setminus B_{r_m})} \le C \| f\|_{L^2(\Omega)}. 
\end{equation}
Let $v_R \in H^1(B_{r_m} \setminus B_{r_1})$ be the reflection of $v$ through $\partial B_{r_m}$, i.e., 
\begin{equation}\label{thm-scheme2-k-est4}
v_R(x) = v \big([r_m - |x|]x/ |x| \big) \mbox{ in } B_{r_m} \setminus B_{r_1}. 
\end{equation}
Define 
\begin{equation*}
u_0  =  \left\{\begin{array}{cl} \hat u &  \mbox{ in } \Omega \setminus B_{r_2}, \\[6pt]
v & \mbox{ in } B_{r_2} \setminus B_{r_m}, \\[6pt]
v_R & \mbox{ in } B_{r_m} \setminus B_{r_1}, \\[6pt]
\hat u (r_2 \cdot / r_1) & \mbox{ in } B_{r_1}. 
\end{array}\right.
\end{equation*}
Then  $u_0 \in H^1\big(\Omega \setminus (\partial B_{r_1} \cup \partial B_{r_2}) \big)$
\begin{equation}\label{part1-u0-2}
\dive (A_0 \nabla u_0) + k^2 \Sigma_0 u_0 = f \mbox{ in } \Omega \setminus (\partial B_{r_1} \cup \partial B_{r_2}). 
\end{equation}
On the other hand, from the definition of $u_0$ and $v$, we have
\begin{equation}\label{part2-u0-2}
[u_0] = [A_0 \nabla u_0 \cdot e_r ] = 0 \mbox{ on } \partial B_{r_2}.
\end{equation}
The properties of the reflection and the definition of $A^H$ garantee
that the transmission conditions also hold on $\partial B_{r_m}$, 
and from the definition of $\hat A$ and \eqref{WE-thm2-bdry}, we obtain 
\begin{equation}\label{part3-u0-2}
[u_0] = [ A_0 \nabla u_0 \cdot e_r ] = 0 \mbox{ on } \partial B_{r_1}. 
\end{equation}
A combination of \eqref{part1-u0-2}, \eqref{part2-u0-2}, and \eqref{part3-u0-2} yields that $u_0 \in H^1(\Omega)$ and satisfies 
\begin{equation*}
\dive (A_0 \nabla u_0) + k^2 \Sigma_0 u_0 = f \mbox{ in } \Omega; 
\end{equation*}
which implies that $u_0$ is a solution of \eqref{sys-k-scheme2} with $\delta = 0$. We also obtain from \eqref{thm-scheme2-k-est0}, \eqref{thm-scheme2-k-est1}, \eqref{thm-scheme2-k-est2}, and \eqref{thm-scheme2-k-est3} that
\begin{equation} \label{thm-scheme2-k-est4}
\| u_0\|_{H^1(\Omega)} \le C \| f\|_{L^2(\Omega)}. 
\end{equation}

We next establish the uniqueness of $u_0$. Let $w_0 \in H^1(\Omega)$ be a solution 
to \eqref{sys-k-scheme2} with $\delta~=~0$. 
Expanding $w$ in spherical harmonics shows
that this function is fully determined in $B_{r_2} \setminus B_{r_m}$
from the Cauchy data 
$w(r_2 \hat{x}), \partial_r w(r_2 \hat{x}), \hat{x} \in \partial B_1$.
Given the form of the coefficients $A^H$, $w$ must also have the symmetry
\begin{equation*}
w_0(x) = w_0 \big([r_m - |x|]x/ |x| \big) \mbox{ in } B_{r_m} \setminus B_{r_1}. 
\end{equation*}
It follows that for $\hat{x} \in \partial B_1$
\begin{eqnarray*}
w_0(r_2 \hat{x}) \;=\; w_0(r_1 \hat{x})
&\textrm{and}&
\partial_r w_0(r_2 \hat{x}) \;=\; \partial_r w_0(r_1 \hat{x}).
\end{eqnarray*}
Thus the function $\hat{w}$ defined by
\begin{eqnarray*}
\hat{w}(x) &=&
\left\{ \begin{array}{ll}
w_0(x) & x \in \Omega \setminus B_{r_2}
\\[6pt]
w_0(r_1 x/r_2) & x \in B_{r_2},
\end{array} \right.
\end{eqnarray*}
is a solution to~\eqref{def-hA-hu-2}.
By uniqueness, $\hat{w} = \hat{u}$, which in turn implies that 
$w_0 = u_0$, which yields the  uniqueness. 
\medskip

Finally, we establish \eqref{thm-scheme2-k-statement1}. Set 
\begin{equation}\label{def-vd2}
v_\delta = u_\delta - u_0 \mbox{ in } \Omega. 
\end{equation} 
It is easy to see that $v_\delta \in H^1_0(\Omega)$ and that it satisfies 
\begin{equation*}
\dive(A_\delta \nabla v_\delta) + k^2 \Sigma_\delta v_\delta =  \dive (i \delta   \mathds{1}_{B_{r_2} \setminus B_{r_1}} \nabla u_0)  - i \delta k^2 \mathds{1}_{B_{r_2} \setminus B_{r_1}} u_0 \mbox{ in } \Omega. 
\end{equation*}
Applying \eqref{lem1-part2} of  Lemma~\ref{lem1}, we derive from \eqref{thm-scheme2-k-est4} that  
\begin{equation}\label{thm-scheme2-k-est5}
 \| v_\delta \|_{H^1(\Omega)} \le C \| \nabla u_0\|_{L^2(\Omega)},
\end{equation}
which is the uniform bound in \eqref{thm-scheme2-k-statement1}. Applying \eqref{lem1-part1} of Lemma~\ref{lem1} and using \eqref{thm-scheme2-k-est4} and \eqref{thm-scheme2-k-est5}, we obtain 
$$
\|u_\delta - u_0\|_{H^1(\Omega)}^2 = \|v_\delta \|_{H^1(\Omega}^2 
\le C \left| \int_{B_{r_2} \setminus B_{r_1}} i \nabla u_0 \nabla v_\delta \right| 
\rightarrow 0, 
$$ 
as $v_\delta$ converges weakly to $0$, which completes the proof. 
\proofend

\begin{remark}  \fontfamily{m} \selectfont  
The proof of Theorem~\ref{thm-main}, where $A^H$ is given by \eqref{scheme2}, 
follows similarly and is left to the reader.
\end{remark}


\section{Constructing hyperbolic metamaterials}

In this section, we show how one can design the type of hyperbolic media 
used in the previous sections, by homogenization of layered materials.
We restrict ourselves to superlensing using HMMs via complementary property 
in the three dimensional quasistatic case, in order to build
a medium $A^H_{\delta}$ that satisfies, as $\delta \to 0$,
\begin{eqnarray} \label{def_AH}
A^H_{\delta} &\rightarrow&
A^H \;=\;
\left\{ \begin{array}{cl} 
\dsp \frac{1}{r^{2}}e_r \otimes e_r 
- (I - e_r \otimes e_r) 
& \textrm{ in}\; B_{r_2} \setminus B_{r_m}, 
\\[6pt]
\dsp -  \frac{1}{r^2} e_r \otimes e_r 
+ (I - e_r \otimes e_r)  
& \textrm{in}\;  B_{r_m} \setminus B_{r_1}, 
\end{array} \right.
\end{eqnarray}
such as that considered in~\eqref{scheme2}. Recall that $r_m  = (r_1 + r_2)/2$. 
The argument can easily be adapted to tuned superlensing using HMMs
in two dimensions  and to  superlensing using HMMs via complementary property in two dimensions 
and to the finite frequency regime.
Our approach follows the arguments developped by 
Murat and Tartar~\cite{KohnCherkaev} for the homogenization of laminated composites.
\medskip

For a fixed  $\delta > 0$, let $\theta = 1/2$ and let $\chi$ denote the characteristic function
of the interval $(0,1/2)$. For $\eps >0$, set, for  $x \in B_{r_2} \setminus B_{r_m}$, 
 \begin{eqnarray*}
b_{1,\eps,\delta}(x) &=& \frac{1}{r^2} 
\left[ (-1 - i\delta) \chi(r/\eps) + \big(1 - \chi(r/\eps) \big)/3 \right] 
\\
b_{2,\eps,\delta}(x) &=&  
 (-3 - i\delta) \chi(r/\eps) + \big(1 - \chi(r/\eps) \big),
\end{eqnarray*}
and, for $x \in B_{r_m} \setminus B_{r_1}$, 
\begin{eqnarray*}
b_{1,\eps,\delta}(x) &=& 
\frac{1}{r^2} \Big[ (-1/3 - i\delta) \chi(r/\eps) + \big(1 - \chi(r/\eps) \big) \Big]
\\
b_{2,\eps,\delta}(x) &=&(-1 - i\delta) \chi(r/\eps) + 3\big(1 - \chi(r/\eps) \big). 
\end{eqnarray*}
Note that since periodic functions converge weakly* to their average in $L^\infty$,
one can easily compute the $L^\infty$ weak-* limits
\begin{equation}\label{def_bH}
b_{1,H,\delta} := \left( w*-\lim_{\eps \to 0} ( b_{1,\eps, \delta})^{-1} \right)^{-1} \quad \mbox{ and } \quad b_{2,H,\delta} :=  w*-\lim_{\eps \to 0} b_{2,\eps, \delta}, 
\end{equation}
and in particular we have in $B_{r_2} \setminus B_{r_m}$
\begin{eqnarray}\label{conv_beps_1}
\left\{\begin{array}{lclcl}
b_{1,H,\delta}(x) &=& \displaystyle\frac{2(1+ i \delta)}{r^2(2 + 3i\delta)}
&=& \left(1 - i\delta/2 + O(\delta^2)\right)/r^2, 
\\[6pt]
b_{2,H,\delta}(x)&=& ( -1 - i \delta/2),
\end{array}\right.
\end{eqnarray}
and in $B_{r_m} \setminus B_{r_1}$
\begin{eqnarray}\label{conv_beps_2}
\left\{\begin{array}{lclcl}
b_{1,H,\delta} (x) &=& \dsp\frac{-2/3 - 2 i\delta}{r^2(2/3 - i\delta)}
&=& -1 - 9i\delta/2 + O(\delta^2), 
\\[6pt]
b_{2,H,\delta} (x)&=& ( 1 - i \delta/2).
\end{array}\right.
\end{eqnarray}
Set 
\begin{equation}\label{def-aeps}
a_{\eps,\delta}(x) =
b_{1,\eps,\delta}(r) e_r \otimes e_r + 
b_{2,\eps,\delta}(r) \left( e_\theta \otimes e_\theta + e_\varphi \otimes e_\varphi \right).  
\end{equation}
Let $a$ be a uniformly elliptic matrix-valued function and 
define 
\begin{eqnarray}\label{Aepsdelta}
A_{\eps,\delta}(x) &=& \left\{ \begin{array}{cl}
I & \mbox{ in } \Omega \setminus B_{r_2}, 
\\[6pt]
a_{\eps,\delta} & \mbox{ in } B_{r_2} \setminus B_{r_1}, 
\\[6pt]
a & \mbox{ in } B_{r_1}, 
\end{array} \right.
\end{eqnarray}
and 
\begin{eqnarray}\label{AHdelta}
A_{\delta}^H(x) &=& \left\{ \begin{array}{cl}
I & \mbox{ in } \Omega \setminus B_{r_2}, 
\\[6pt]
b_{1,H,\delta}  e_r \otimes \, e_r + b_{2,H,\delta} \, \left( e_\theta \otimes e_\theta + e_\varphi \otimes e_\varphi \right) & \mbox{ in } B_{r_2} \setminus B_{r_1}, 
\\[6pt]
a & \mbox{ in } B_{r_1}.
\end{array} \right.
\end{eqnarray}

We have 

\begin{proposition}~\label{prop_constr}
Let $0 < r_1 < r_2$, and let $\Omega$ be a smooth bounded connected open subset of $\mR^3$ 
such that 
$B_{r_2} \subset \subset \Omega$. Given  $f \in L^2(\Omega)$ with $\supp f \cap B_{r_2}  = \O$,  let $u_{\e,\delta} \in H^1_0(\Omega)$ be the unique solution to
\[
\dive (A_{\eps,\delta} \nabla u_{\eps, \delta}) = f  \mbox{ in }  \Omega,
\]
where $A_{\eps, \delta}$ is given by \eqref{Aepsdelta}. 
Then, as $\eps \to 0$, 
$u_{\eps,\delta}$ converges weakly in $H^1(\Omega)$ to  $u_{H,\delta} \in H^1_0(\Omega)$ the unique solution of the equation 
\[
\dive(A^H_{\delta} u_{H,\delta}) = f \mbox{ in } \Omega,
\]
where $A^H_\delta$ is defined by \eqref{AHdelta}. 
\end{proposition}

\begin{remark}  \fontfamily{m} \selectfont   Materials given in \eqref{def-aeps} could in principle be fabricated as a laminated composite containing anisotropic metallic phases with a conductivity described by a Drude model.
Also note that the imaginary part of $A_\delta^H$ has the form
$-i \delta M$, where $M$ is a diagonal, positive definite matrix,
and is not strictly equal to $- i \delta I$ as in the hypotheses 
of Theorem~\ref{thm-main}.
Nevertheless, its results hold for this case as well.
\end{remark}

\noindent{\bf Proof.} 
For notational ease, we drop the dependance on $\delta$ in the notation.
By Lemma~\ref{lem1} (see also Remark~\ref{rem-D-0}),  there exists a unique solution 
$u_{\eps} \in H^1_0(\Omega)$ to
\begin{equation}\label{eq-ueps}
\textrm{div}(A_\eps \nabla u_\eps) = f  \mbox{ in } \Omega, 
\end{equation}
which further satisfies $||u_\eps||_{H^1(\Omega)} \;\leq\; C\, ||f||_{L^2(\Omega)}$, 
with $C$ independent of $\eps$ (it may depend on $\delta$ though).
We may thus assume, that up to a subsequence, $u_\eps $ converges weakly in $H^1(\Omega)$
to some $u_H \in H^1(\Omega)$.
Standard results in homogenization~\cite{KohnCherkaev} show that $u_H \in H^1_0(\Omega)$ solves
an equation of the same type as~\eqref{eq-ueps}: 
\begin{equation}\label{eq_ueps2}
\dive (A^H \nabla u_H) = f  \mbox{ in } \Omega, 
\end{equation}
where the tensor of homogenized coefficients $A_H$ has the form
\begin{eqnarray*}
A_H(x) &=& \left\{ \begin{array}{cl}
I & \mbox{ for } x \in \Omega \setminus B_{r_2}, 
\\[6pt]
a_H(x) &  \mbox{ for } x \in B_{r_2} \setminus B_{r_1}, 
\\[6pt]
a(x) &  \mbox{ for } x \in B_{r_1}.
\end{array} \right.
\end{eqnarray*}
To identify the tensor $a_H$, set
\begin{equation} \label{def_sigma}
\sigma_{1,\eps} = r^2 b_{1,\eps} \partial_r u_\eps
\quad \textrm{in}\; B_{r_2} \setminus B_{r_1}. 
\end{equation}
Using spherical coordinates in $B_{r_2} \setminus B_{r_1}$, we have 
$$
\dive (A_\eps \nabla u_\eps) = \frac{1}{r^2} \partial_r (r^2 b_{1, \eps} \partial_r u_\eps) + \frac{b_{2, \eps}}{r^2} \Delta_{\partial B_{1}} u_\eps  \mbox{ in } B_{r_2} \setminus B_{r_1}, 
$$
where $\Delta_{\partial B_{1}}$ denotes the Laplace-Beltrami operator on $\partial B_1$. 
This implies, since $\supp f \cap B_{r_2} = \O$,  
\begin{equation*}
\partial_r \sigma_{1,\eps} =  - \Delta_{\partial B_1} \big( b_{2,\eps}(r) u_{\eps} \big)
 \mbox{ in } B_{r_2} \setminus B_{r_1},
\end{equation*}
since $b_{2,\eps}$ only depends on $r$ for a fixed $\eps$.
Consequently, 
$\sigma_{1,\eps}$ and $\partial_r \sigma_{1,\eps}$ are uniformly bounded
with respect to $\eps$ in 
$L^2\big(r_1, r_2, L^2(\partial B_1) \big)$ and in
$L^2\big(r_1, r_2, H^{-1}(\partial B_1) \big)$ respectively.
Invoking Aubin compactness theorem as in~\cite{KohnCherkaev}, we infer that up to a subsequence, 
$\sigma_{1,\eps}$ converges strongly in 
$L^2\big(r_1, r_2, H^{-1}(\partial B_1) \big)$
to some limit $\sigma_{1,H} \in L^2(B_{r_2} \setminus B_{r_1})$.
Rewriting~\eqref{def_sigma} as
\begin{equation*} 
\left( r^2 b_{1,\eps} \right)^{-1} \sigma_{1,\eps} =  \partial_r u_\eps,
\end{equation*}
and letting  $\eps \to 0$, yields
\begin{eqnarray*}
\sigma_{1,H} &=&
\left( w*-\lim (r^2 b_{1,\eps})^{-1} \right)^{-1} \partial_r u_H
\\
&=&
\displaystyle \frac{r^2}{ w*-\lim ( b_{1,\eps})^{-1} }  \partial_r u_H.
\end{eqnarray*}
On the other hand, since $u_\eps \to u_H$ strongly in $L^2(\Omega)$, it 
follows that $b_{2,\eps}(r) u_\eps \to w*-\lim b_{2,\eps}(r) u_H$
in $L^2$. We derive that
\begin{equation} \label{eq_hom}
\partial_r \big( r^2 b_{1,H} \partial_r u_H \big)
+ \Delta_{\partial B_1} \big(b_{2,H} u_H \big)
=0
\mbox{ in } B_{r_2} \setminus B_{r_1},
\end{equation}
where
$b_{1,H} \;=\; \left( w*-\lim ( b_{1,\eps})^{-1} \right)^{-1}$
and $b_{2,H} \;=\; w*-\lim b_{2,\eps}$.
We can then identify
\begin{equation*}
a_H = b_{1,H} e_r \otimes e_r 
+ b_{2,H} \left( e_\theta \otimes e_\theta + e_\varphi \otimes e_\varphi \right),
\end{equation*}
which, given~(\ref{conv_beps_1}--\ref{conv_beps_2}), has the form considered in~\eqref{scheme2-k}. 
\medskip

Since periodic functions weakly-* converge to their average in $L^\infty$
one easily checks that in fact the whole sequence $u_\eps$ converges to the unique $H^1_0$-solution to~\eqref{eq_hom}. \proofend


\providecommand{\bysame}{\leavevmode\hbox to3em{\hrulefill}\thinspace}
\providecommand{\MR}{\relax\ifhmode\unskip\space\fi MR }
\providecommand{\MRhref}[2]{%
  \href{http://www.ams.org/mathscinet-getitem?mr=#1}{#2}
}
\providecommand{\href}[2]{#2}

\end{document}